\documentstyle[11pt]{article}
\newtheorem{theorem}{Theorem}
\newtheorem{definition}[theorem]{Definition}

\newtheorem{example}[theorem]{Example}

\newtheorem{proposition}[theorem]{Proposition}

\newtheorem{odstavec}{\hskip -2mm}[section]

\catcode`\@=11

\def\ps@myheadings{\let\@mkboth\@gobbletwo
\def\@oddhead{\ifnum\count0=1 \hfill\else
\rightmark \hfil \rm\thepage\fi}%
\def\@oddfoot{\ifnum\count0=1 \hfill \rm 1 \hfill \else
\hfill\fi}
\def\@evenhead%
{\rm\leftmark\hfil\rm\thepage}%
\def\@evenfoot{}\def\sectionmark##1{}
\def\subsectionmark##1{}}

\def\@begintheorem#1#2{\it \trivlist \item[\hskip
 \labelsep{\bf #1\ #2.}]}
\def\@opargbegintheorem#1#2#3{\it \trivlist\item[\hskip%
 \labelsep{\bf #1\ #2.\ (#3)}]}
\def\@endtheorem{\endtrivlist}

\def\@listI{\leftmargin\leftmargini \parsep 1pt plus 2.5pt
 minus 1pt\topsep 5pt plus 2pt minus 3pt%
 \itemsep 0pt plus 2.5pt minus 1pt}
\let\@listi\@listI
\@listi

\def\@sect#1#2#3#4#5#6[#7]#8{\ifnum #2>\c@secnumdepth%
 \def \@svsec {}\else \refstepcounter {#1}\edef \@svsec%
 {\csname the#1\endcsname. \hskip .1em }\fi \@tempskipa%
 #5\relax \ifdim \@tempskipa >\z@ \begingroup #6\relax%
 \@hangfrom {\hskip #3\relax \@svsec }{\interlinepenalty%
 \@M #8.\par }\endgroup \csname #1mark\endcsname {#7}%
 \addcontentsline {toc}{#1}{\ifnum #2>\c@secnumdepth%
 \else \protect \numberline {\csname the#1\endcsname. }%
 \fi #7}\else \def \@svsechd {#6\hskip #3\@svsec #8.%
 \csname #1mark\endcsname {#7}\addcontentsline {toc}{#1}%
 {\ifnum #2>\c@secnumdepth \else \protect \numberline%
 {\csname the#1\endcsname. }\fi #7}}\fi \@xsect {#5}}

\def\section{\@startsection {section}{1}{\z@ }%
 {-3.5ex plus -1ex minus -.2ex}{2.3ex plus .2ex}{\bf }}

\def\thebibliography#1{%
 \section *{\hfil REFERENCES \@mkboth {REFERENCES}{REFERENCES}}%
 \list {[\arabic {enumi}]}{\settowidth \labelwidth {[#1]}%
 \leftmargin \labelwidth \advance \leftmargin \labelsep %
 \usecounter {enumi}} \def \newblock %
 {\hskip .11em plus .33em minus -.07em} \sloppy \clubpenalty 4000%
 \widowpenalty 4000 \sfcode`\.=1000\relax}

\def\@maketitle{%
 \newpage \null \vskip 2em
 \begin{center}
{\Large\bf \@title \par }
 \vskip 1.5em
 {\large \lineskip .5em
 \begin {tabular}[t]{c}\@author
 \end{tabular}\par}
 \end{center}
  \vskip .8em}

\def\abstract{%
\if@twocolumn \section *{Abstract}
 \else \small\quotation\noindent{\bf Abstract.}\fi}

\catcode`\@=13

\font\tenbold=msbm10 scaled \magstep1
\font\sevenbold=msbm7 scaled \magstep1
\font\fivebold=msbm5 scaled \magstep1
\newfam\boldfam \def\Bbb{\fam\boldfam\tenbold}
\textfont\boldfam=\tenbold
\scriptfont\boldfam=\sevenbold
\scriptscriptfont\boldfam=\fivebold

\font\tengoth=eufm10 scaled \magstep1
\font\sevengoth=eufm7 scaled \magstep1
\font\fivegoth=eufm5 scaled \magstep1
\newfam\gothfam \def\frak{\fam\gothfam\tengoth}
\textfont\gothfam=\tengoth
\scriptfont\gothfam=\sevengoth
\scriptscriptfont\gothfam=\fivegoth

\font\tenmale=cmti10
\font\sevenmale=cmti7
\font\fivemale=cmti5 
\newfam\malefam \def\male{\fam\malefam\tenmale}
\textfont\malefam=\tenmale
\scriptfont\malefam=\sevenmale
\scriptscriptfont\malefam=\fivemale

\def\pa{{\partial}}

\def\calF{{\cal F}}
\def\calP{{\cal P}}
\def\Riso{{\cal R}_{\rm iso}}
\def\alphaiso{{\alpha_{\rm iso}}}

\def\pretRiso{{\it pre}{\tilde {\cal R}_{\rm iso}}}
\def\tRiso{{\tilde {\cal R}_{\rm iso}}}
\def\Dif{{\cal D}{\it if}}
\def\preDif{\mbox{\it pre}{\cal D}{\it if}}

\def\F{{\male F}}
\def\H{{\male H}}
\def\G{{\male G}}
\def\idmala{{\male 1\hskip -1.8mm 1}} 
\def\L{{\male L}}

\def\FiltOp{{\tt FiltOp_{{\Bbb Z}}}}
\def\vlra{{\hbox{$-\hskip-1mm-\hskip-2mm\longrightarrow$}}}
\def\Eout{{E_{\rm out}}}
\def\SHE{({\underline \F,\underline \G,\underline \H,\underline \L})}
\def\tSHE{({\tilde {\underline \F},\tilde {\underline \G},
            \tilde {\underline \H},\tilde {\underline \L}})}
\def\rfake{{r_{\rm fake}}}
\def\otexp#1#2{{#1^{\otimes #2}}}
\def\iotafake{{\iota_{\rm fake}}}
\def\tF{{\tilde \F}}
\def\tG{{\tilde \G}}
\def\tH{{\tilde \H}}
\def\tL{{\tilde \L}}
\def\tcalRfake{{\tilde {\cal R}_{\rm fake}}}
\def\Ain{{A_{\rm in}}}
\def\Aout{{A_{\rm out}}}
\def\td{{\tilde d}}
\def\calRfake{{{\cal R}_{\rm fake}}}
\def\kernel{{\frak Z}}
\def\udeg{{\underline{\it de}\hskip -.5mm {\it g}}}
\def\fbar{{\overline f}}
\def\gbar{{\overline g}}
\def\catChain{{\tt Chain_{{\Bbb Z}}}}

\def\colorop #1(#2;#3){{#1}
   \left(\rule{0pt}{15pt}\right.
         \hskip -3mm \begin{array}{c}
	              #3\\#2
                     \end{array}
         \hskip -3mm \left. 
   \rule{0pt}{15pt} \right)
}
\def\doubless#1#2{{
\def\arraystretch{.5}
\begin{array}{c}
\mbox{\scriptsize $\scriptstyle #1$}
\\
\mbox{\scriptsize $\scriptstyle #2$}
\end{array}\def\arraystretch{1}
}}
\def\FF#1{{\cal F}(f_{< #1}; g_{< #1})}
\def\cases#1#2#3#4{
                  \left\{
                         \begin{array}{ll}
                           #1,\ &\mbox{#2}
                           \\
                           #3,\ &\mbox{#4}
                          \end{array}
                   \right.
}
\def\Hom{{\it Hom}}
\def\frakO{{\frak o}}
\def\Ares{{A}_{\rm res}}

\def\alphafake{{\alpha}_{\rm fake}}
\def\B{{\tt B}}
\def\W{{\tt W}}
\def\kernelb{{\kernel_\bullet}}
\def\qed{\hspace*{\fill}
         \mbox{\hphantom{mm}\rule{0.25cm}{0.25cm}}\\}
\def\xb{{\overline x}}
\def\yb{{\overline y}}
\def\id{1\!\!1}

\def\Iso{{\cal I}{\it so}}

\def\deg{{\it deg}}
\def\End{\hbox{${\cal E}\hskip -.1em {\it nd}$}}
\def\seq#1#2{\{#1\}_{{#2}}}

\def\widedef{{\hskip 4mm := \hskip 4mm}}

\def\freeprod{{\Dif*\Riso}}
\def\hb{h_\bullet}
\def\lb{l_\bullet}
\def\fb{f_\bullet}
\def\gb{g_\bullet}
\def\hbtilde{{\tilde h}_\bullet}
\def\lbtilde{{\tilde l}_\bullet}
\def\fbtilde{{\tilde f}_\bullet}
\def\gbtilde{{\tilde g}_\bullet}

\begin{document}
\baselineskip15pt
\overfullrule4pt

                                                      %
\newpage \null \vskip 2em\vskip.5cm                   %
\begin{center}{\Large\bf IDEAL PERTURBATION LEMMA     %
 \par }                                               %
\vskip 1.5em                                          %
{\large \lineskip .5em                                %
Martin Markl}\footnote{Supported by the               %
grant GA \v CR 201/99/0675}                           %
\vskip .8em                                           %
Mathematical Institute of the Academy, \v Zitn\'a 25  %
\\                                                    %
115 67 Praha 1, The Czech Republic,                   %
\catcode`\@=11                                        %
\\                                                    %
email: {\tt markl@math.cas.cz}                        %
\catcode`\@=13                                        %
\end{center}\par \vskip 1.5em                         %
                                                      %

\begin{abstract}
We explain the essence of perturbation
problems. The key to understanding is the structure
of chain homotopy equivalence -- the standard one
must be replaced by a finer notion which we call a {\em strong\/}
chain homotopy equivalence.

We formulate an Ideal Perturbation Lemma 
and show how both new and classical (including the Basic Perturbation
Lemma)  results follow from this
ideal statement.
\end{abstract}

\section{Introduction and results}
\label{Introduction}

All algebraic objects in this paper are defined over the ring of
integers ${\Bbb Z}$.
Our work was motivated by reconsidering the following classical Basic
Perturbation Lemma (see~\cite{huebschmann-kadeishvili:MZ91} 
and the historical account there). 

Recall that
{\em strong deformation retract
(SDR)\/} data (also called a {\em contraction\/}) are given by 
chain complexes $(M,d_M)$, $(N,d_N)$, 
chain maps 
$\F : (M,d_M) \to (N,d_N)$, $\G : (N,d_N) \to (M,d_M)$ 
and a chain homotopy $\H : M
\to M$ satisfying
\begin{eqnarray}
\label{SDR}
\F d_M &=& d_N \F,\nonumber \\
\G d_N &=& d_M \G,\\
\G\F - \idmala_M &=& d_M \H + \H d_M \mbox { and }\nonumber \\
\nonumber \F\G &=& \idmala_N.
\end{eqnarray}

This of course means that $(N,d_N)$ is a strong deformation retract of
$(M,d_M)$.
One usually
assumes that the following {\em side conditions\/} (also called {\em
annihilation properties\/}) hold:
\begin{equation}
\label{side_conditions}
\H \H=0,\ \H \G = 0 \mbox { and } \F \H =0.
\end{equation}
Then the following statement is true.

\noindent 
{\bf Basic Perturbation Lemma (BPL).}{\em\
Suppose we are given strong deformation data~(\ref{SDR})
satisfying~(\ref{side_conditions}) and a
perturbation ${\tilde d}_M$ of the differential $d_M$
on M. Then there are 
perturbations ${\tilde d}_N, {\tilde \F}, {\tilde \G}$  
and $\tilde \H$ of $d_N, \F,\G$ and $\H$ that again form 
strong deformation data~(\ref{SDR}),
\[
\tilde \F \td_M = \td_N \tilde \F, \hskip 3mm
\tilde \G \td_N = \td_M \tilde \G, \hskip 3mm
\tilde \G\tilde \F - \idmala_M = \td_M \tilde \H + \tilde \H \td_M \
\mbox { and }\
\tilde \F\tilde \G = \idmala_N. 
\]
}

All notions used in the formulation of the BPL are standard and we believe
it is not necessary to repeat their definitions here.
Filtered objects and perturbations are treated in
Section~\ref{Filtrations}.
The perturbation $({\tilde d}_N, {\tilde \F}, {\tilde \G},\tilde \H)$
is given by the following explicit formulas 
(see again~\cite{huebschmann-kadeishvili:MZ91}): 
\[
\def\arraystretch{1.4}
\begin{array}{lcr}
{\tilde d}_N &=&
d_N + \F(\pa_M + \pa_M \H \pa_M + \pa_M \H \pa_M \H \pa_M + 
        \pa_M \H \pa_M \H \pa_M \H \pa_M+\cdots )\G,
\\
\tilde \F &=& \F +\F(\pa_M + \pa_M \H \pa_M + \pa_M \H \pa_M \H \pa_M + 
        \pa_M \H \pa_M \H \pa_M \H \pa_M+ \cdots )\H,
\\
\tilde \G &=& \G +\H(\pa_M + \pa_M \H \pa_M + \pa_M \H \pa_M \H \pa_M + 
        \pa_M \H \pa_M \H \pa_M \H \pa_M+ \cdots )\G,
\\
\tilde \H &=& \H +\H(\pa_M + \pa_M \H \pa_M + \pa_M \H \pa_M \H \pa_M + 
        \pa_M \H \pa_M \H \pa_M \H \pa_M+ \cdots )\H,
\end{array}
\]
where $\partial_M := \tilde d_M - d_M$.
The formulas above contain infinite series, so one must assume some
conditions assuring that they converge. This is usually achieved by
assuming that both $(M,d_M)$ and $(N,d_N)$ are filtered complete, see
again Section~\ref{Filtrations}.

Our original motivation was to understand why there is such a formula
and what is the role of side conditions. As usual, the best way to
understand a problem is to formulate it in as general a form as
possible. So let us consider the following:

\vskip 2mm
\noindent 
{\bf Perturbation Problem (PP).}
{\em  Suppose we are given two complete
filtered complexes $M = (N,d_M)$ and $N
= (N,d_N)$ and chain maps $\F : M \to N$ and $\G : N \to M$
that are chain homotopy inverse to each other, with homotopies $\H : M
\to M$ and $\L : N \to N$, that is
\begin{equation}
\label{muj_Misa}
\F d_M \!=\! d_N \F,\
\G d_N \!=\! d_M \G,\
\G\F - \idmala_M \!=\! d_M \H + \H d_M,\
\F\G -\idmala_N \!=\! d_N \L + \L d_N.
\end{equation}
Given a perturbation ${\tilde d}_M$ of the differential $d_M$, find
perturbations ${\tilde d}_N, {\tilde \F}, {\tilde \G}, {\tilde \H}$  
and $\tilde \L$ of $d_N, \F,\G,\H$ and $\L$ such that ${\tilde \F}$ and 
${\tilde \G}$ are chain maps with respect to the perturbed
differentials, 
homotopy inverse to each other, 
with homotopies ${\tilde \H}$ and ${\tilde \L}$, that is
\[
{\tilde \F}{\tilde d}_M \!=\! {\tilde d}_N {\tilde \F},\
{\tilde \G} {\tilde d}_N \!=\!  {\tilde d}_M {\tilde \G},\
{\tilde \G}{\tilde \F}- \idmala_M \!=\! {\tilde d}_M {\tilde \H} + {\tilde \H} {\tilde d}_M,\
{\tilde \F}{\tilde \G} -\idmala_N \!=\! {\tilde d}_N {\tilde \L} + {\tilde \L} {\tilde d}_N.
\]
\/}

Observe that in the formulation of the BPL and the PP we consider not
only the differentials and 
the chain maps, but also the homotopies to be a part of the
structure which has to be perturbed. Ignoring homotopies
leads to the `crude' perturbation lemma formulated at the end of this
Introduction.

The fact, both frustrating and provoking, is that the PP has, for
general input data,
\underline{no} \underline{solution}! -- a
rigorous formulation and proof of this negative statement is provided by 
Theorem~\ref{Harvey}.
The reason is that the chain homotopy equivalence
$(\F,\G,\H,\L)$ of~(\ref{muj_Misa}) 
is not a homotopy invariant concept and it must be replaced by
a subtler notion of a {\em \underline{stron}g 
(chain) homotopy equivalence\/}:

\begin{definition}
\label{Kiki}
A strong homotopy equivalence (SHE) consists of degree $2m$ maps
$\F_{2m} : M \to N$, $\G_{2m} : N \to M$ and degree $2m+1$ `homotopies'
$\H_{2m+1} : M \to M$ , $\L_{2m+1}: N \to N$, for all
$m \geq 0$, such that
\[
\def\arraystretch{1.4}
\begin{array}{rclrcl}
\F_0 d_M& = & d_N \F_0, &
\G_0 d_N& = &d_M \G_0,
\\
\G_0\F_0 - \idmala_M& =& d_M \H_1 + \H_1 d_M, &
\F_0\G_0 -\idmala_N & =& d_N \L_0 + \L_0 d_N
\end{array}
\]
and that, for each $m \geq 1$,
\begin{eqnarray}
\nonumber
d_N \F_{2m} - \F_{2m} d_N   &=& \sum_{0 \leq i < m}(\F_{2i}\H_{2(m-i)-1}-
\L_{2(m-i)-1}\F_{2i}),
\\
\nonumber 
d_M\H_{2m+1}+ \H_{2m+1}d_M &=& \sum_{0 \leq j \leq m} \G_{2j}\F_{2(m-j)} -
               \sum_{0 \leq j < m} \H_{2j+1}\H_{2(m-j)-1},
\\
\nonumber 
d_M \G_{2m} - \G_{2m} d_N   &=& \sum_{0 \leq i < m}(\G_{2i}\L_{2(m-i)-1}-
\H_{2(m-i)-1}\G_{2i}),
\\
\nonumber 
d_N \L_{2m+1} + d_N \L_{2m+1} &=& \sum_{0 \leq j \leq m} \F_{2j}\G_{2(m-j)} -
               \sum_{0 \leq j < m} \L_{2j+1}\L_{2(m-j)-1}.
\end{eqnarray}
\end{definition}

See~\ref{8.1} where we expanded the axioms above for some small $m$. 
To understand
better the meaning of a SHE, we
offer the following analogy.

A {\em homotopy associative algebra\/} is a chain complex $V = (V,d_V)$
with a homotopy associative multiplication $\mu :\otexp V2 \to V$:
\[
\mu(\mu \otimes \id) - \mu(\id \otimes \mu) \sim 0 \mbox { modulo a
chain homotopy $\nu : \otexp V3 \to V$}.
\]
As argued in~\cite{markl:ha},
a proper homotopy invariant version of this concept is that of a
{\em strongly\/} 
homotopy associative algebra, which
is a structure consisting of infinitely many
multilinear operations $\{\mu_n : \otexp Vn \to V\}_{n \geq 2}$
such that the `multiplication' $\mu_2 : \otexp V2 \to V$ is homotopy
associative up to the homotopy $\mu_3 : \otexp V3 \to V$, and there
is, for each $n \geq 4$, a certain `coherence relation' assumed to be
zero modulo the homotopy $\mu_n$, see~\cite{stasheff:TAMS63}. 
While each strongly
homotopy associative algebra defines, by $\mu := \mu_2$ and $\nu :=
\mu_3$, a homotopy associative one, the converse is not true; there
are obstructions for extending a homotopy associative
multiplication to a strongly homotopy associative one.

The situation in Definition~\ref{Kiki} is similar. 
While a strong homotopy equivalence 
defines, by $\F := \F_0$, $\G := \G_0$, $\H := \H_1$
and $\L := \L_1$ an ordinary homotopy equivalence, 
the converse is not true --
there is a primary
obstruction $[{\frak o}]$ for extending a homotopy equivalence to a
strong one. The surprising Theorem~\ref{snad_nejsem_nastydly} 
says that vanishing of this
primary obstruction already implies the existence of the extension.

A strong homotopy equivalence of $M$ and $N$ will be denoted as $\SHE
: M \to N$. Let us formulate our Ideal Perturbation Lemma.

\vskip 2mm
\noindent 
{\bf Ideal Perturbation Lemma (IPL).}
{\em  Suppose we are given two complete filtered complexes $M = (N,d_M)$ and $N
= (N,d_N)$ and a strong homotopy equivalence $\SHE : M \to N$.

Given a perturbation ${\tilde d}_M$ of the differential $d_M$, 
there exist a 
perturbation ${\tilde d}_N$ of the differential $d_N$ and a
perturbation  $\tSHE$  of $\SHE$ which is a strong homotopy
equivalence of the perturbed complexes $(M,\td_M)$ and $(N,\td_N)$. Moreover,
the perturbations $\td_M$ and $\tSHE$ depend functorially on ${\tilde d}_M$
and $\SHE$. 
\/}

The IPL is proved in Section~\ref{Proofs}, see also~\ref{8.6} for
explicit formulas.
As most ideal things, the Ideal Perturbation Lemma is
almost useless. In practice, the input data are formulated only in
terms of an ordinary homotopy equivalence, and the answer is also
expected to be a perturbation of this ordinary homotopy
equivalence. 
Here is our mundane version of the Ideal Perturbation Lemma.

\begin{theorem}
\label{Micinka}
Suppose that the obstruction $[{\frak o}]$ to the extension of the
homotopy equivalence~(\ref{muj_Misa}) to a strong one
vanishes. 
Then the Perturbation Problem has a solution, functorial up to
a choice of the extension of~(\ref{muj_Misa}) to a strong
homotopy equivalence.
\end{theorem}

The theorem immediately follows from the IPL and the above notes.
There are situations when the obstruction $[{\frak o}]$ vanishes
and when there even exists a functorial extension of the homotopy
equivalence~(\ref{muj_Misa}) to a SHE.
This the case of our motivating example of the Basic Perturbation Lemma
(the case $\L=0$). It immediately follows from 
Theorem~\ref{Marcelka} that the side
conditions~(\ref{side_conditions}) 
guarantee the existence of a functorial extension
of~(\ref{SDR})
to a strong homotopy equivalence. So Theorem~\ref{Micinka} implies the BPL.

Another trick that overrides the nonexistence of a solution of the PP
is to change the initial data a bit. We show in Theorem~\ref{3+1} 
that changing in~(\ref{muj_Misa}) the homotopy $\H$ to
$\H-\G(\F\H-\L\F)$ (or, dually, $\L$ to $\L-\F(\G\L-\H\G)$)
annihilates the obstruction $[\frak o]$ 
and we reprove the following recent result by J.~Huebschmann and
T.~Kadeishvili~\cite{huebschmann-kadeishvili:MZ91}.

\begin{theorem}
\label{zitra}
Let $M = (M,d_M)$ and $N = (N,d_N)$ be complete filtered chain complexes and
$(\F,\G,\H,\L)$ a chain homotopy equivalence~(\ref{muj_Misa}). 

Given a perturbation $\td_M$ of $d_M$, there exist a perturbation
$\td_N$ of the differential $d_N$ and a homotopy equivalence
$(\tF,\tG,\tH,\tL)$ of the perturbed complexes $(M,\td_M)$ and
$(N,\td_N)$ that is a perturbation of $(\F,\G,\H-\G(\F\H-\L\F),\L)$.
\end{theorem}

Changing $\L$ to $\L-\F(\G\L-\H\G)$ and leaving $\H$ untouched gives
the following complement to
Theorem~\ref{zitra}.

\vskip2mm
\noindent 
{\bf Complement to Theorem~\ref{zitra}.}
{\em  
Under the assumption of Theorem~\ref{zitra}, there exists another
perturbation 
$\td_N'$ of the differential $d_N$ and another homotopy equivalence
$(\tF',\tG',\tH',\tL')$ of the perturbed complexes $(M,\td_M)$ and
$(N,\td_N')$ that is a perturbation of $(\F,\G,\H,\L-\F(\G\L-\H\G))$.
\/}

Ignoring the homotopies in the Perturbation Problem, we get the following

\vskip 2mm
\noindent
{\bf Crude Perturbation Lemma.}
{\em  Suppose we are given two complete
filtered complexes $M = (N,d_M)$ and $N
= (N,d_N)$ and chain maps $\F : M \to N$ and $\G : N \to M$
that are chain homotopy inverse to each other.

Given a perturbation ${\tilde d}_M$ of the differential $d_M$, there are
perturbations ${\tilde d}_N, {\tilde \F}$ and ${\tilde \G}$
of $d_N, \F$ and $\G$ such that ${\tilde \F}$ and 
${\tilde \G}$ are chain maps with respect to the perturbed
differentials, 
homotopy inverse to each other.
\/}

A conceptual explanation of these results is given in
Section~\ref{conceptual}.

\begin{center}
-- -- -- -- --
\end{center}

\noindent 
{\bf Plan of the paper:}
In Section~\ref{Language} we recall colored operads
and introduce the operad $\Iso$ describing isomorphisms of chain
complexes. 
In Section~\ref{Filtrations} we repeat necessary facts
on filtrations and perturbations and define the filtered operad
$\Dif$ describing perturbations of differentials.
The filtered operad $\Riso$ that
describes strong homotopy equivalences is introduced in 
Section~\ref{Strong} where we also discuss extensions of a
homotopy equivalence to a strong one. 
In Section~\ref{The_retraction} we introduce the
operad $\tRiso$ for perturbations of strong homotopy
equivalences and construct a retraction $r$ that gives the
functorial solution to the IPL. 
Some of the proofs are postponed to
Section~\ref{Proofs}.
In Section~\ref{conceptual} we give a conceptual
explanation of the results. %
In the Appendix (Section~\ref{Appendix}) we present some explicit formulas.

\section{Language of operads}
\label{Language}

Roughly speaking, operads are objects that describe types of algebraic
systems. Colored operads are then objects describing diagrams of
algebraic systems. The definition of a (colored) operad is
classical (see~\cite{boardman-vogt:73} or~\cite{markl:ha})
and we will not repeat it here in its full
generality.

By an operad we will always mean an operad
in the symmetric monoidal category $\catChain$
of differential graded complexes of
abelian groups (that is, complexes of ${\Bbb Z}$-modules).
Operads in this category
behave in many aspects as associative algebras, so we may
speak about suboperads, ideals, presentations, resolutions, 
etc., see~\cite{ginzburg-kapranov:DMJ94}. 

All algebraic objects in this paper will have only unary
operations. Colored operads describing algebraic systems with only unary
operations are the same as
small additive categories enriched over $\catChain$. This means that all
hom-sets are chain complexes and composition maps are homomorphisms of
chain complexes. All operads in this paper will be of this type.

We will use the `operadic' notation and terminology.
Thus, for such an operad/category $\calP$, we call ${\frak
C} : = {\it Ob}(\calP)$ the {\em set of colors\/} and, for $c,d \in
{\frak C}$, we denote 
\[
\colorop{\calP}(c;d) := 
{\it Mor}_{\calP}(c,d). 
\]
We will usually express the fact that $
f \in \colorop{\calP}(c;d)$ by writing $f : c \to d$. 

In the particular case when ${\rm card}({{\frak C}}) = 1$, the ${\frak
C}$-colored operads are exactly differential graded associative unital
algebras. In this paper, by a
colored operad we always mean an operad colored by 
the two-point set ${\frak C} =
\{\B,\W\}$ ($\B$ from black, $\W$ from white) or by a set isomorphic
to this one.

\begin{example}
\label{Jukl}{\rm
Let $M = (M,d_M)$ be a chain complex, then the {\em endomorphism
operad\/} $\End_M$ is defined to be the chain complex $\Hom(M,M)$ with
the operadic structure (which in this particular case is the same as
that of an unital associative algebra)
given by the composition. An {\em algebra
over an operad\/} $\calP$ is an operadic 
homomorphism $A : \calP \to
\End_M$. In this situation we also say that the operad
$\calP$ {\em acts\/} on the chain complex $M$.
}\end{example}

\begin{example}
{\rm
This example describes a colored version of the endomorphism operad
recalled in Example~\ref{Jukl}.
Let $M = (M,d_M)$ and $N = (N,d_N)$ be chain complexes.
By a 
{\em colored endomorphism operad\/} $\End_{M,N}$ we mean the full
subcategory of $\catChain$ with objects $M$ and $N$. If $\calP$ is a
$\{\B,\W\}$-colored operad, then by a $\calP$-algebra we mean a
homomorphism $A : \calP \to \End_{M,N}$ such that $A(\B) = M$ and
$A(\W) = N$.
}
\end{example}

\begin{example}
\label{trochu}
{\rm
Let $f: \B \to \W$, $g: \W \to \B$ be two degree-zero generators and denote 
\[
\Iso := \left(\frac{\calF(f,g)}
                  {(fg=1_\W,gf=1_\B)}, d=0\right).
\]
In the above display, $\calF(f,g)$ denotes the free 
$\{\B,\W\}$-colored operad on the set
$\{f,g\}$ and $(fg=1_\W,gf=1_\B)$ the operadic ideal
generated by $fg-1_\W$ and $gf-1_\B$.

An algebra $A : \Iso \to \End_{M,N}$
consists of two degree zero chain maps $\F : M \to N$, $\G : N \to M$ that
are inverse to each other. 
Thus the operad $\Iso$ describes isomorphisms of chain complexes,
whence its name.
}
\end{example}

\section{Filtrations and Perturbations}
\label{Filtrations}

Let $M = (M,d_M)$ be a chain complex.
A (descending) {\em filtration\/} on $M$ is a descending
sequence $\{F^p M\}_{p \geq 0}$ of subcomplexes of $M$.
If not stated otherwise, we always assume that the filtration is {\em
complete\/}. This, by definition, means that the module $M$ is
complete in the $F^p$-adic topology. 
This guarantees that each sum $\sum_{p \geq 0} m_p$ with $m_p \in
F^pM$ represents a unique element of $M$. A typical example 
is the module of power series ${\Bbb Z}[[h]]$
with the filtration defined by $F^p {\Bbb Z}[[h]] := h^p {\Bbb
Z}[[h]]$, $p \geq 0$.

Morphisms of filtered chain complexes are maps that preserve 
filtrations. A linear map
$g: M \to N$ is a {\em perturbation\/} or {\em deformation\/} of a
linear map $f : M \to N$ if 
\[
(f - g)(F^pM) \subset F^{p+1}N
\mbox { for each $p \geq 0$.}
\]
If $M$ and $N$ are  filtered complexes, then the chain complex ${\it
Hom}(M,N)$ is also filtered, by
\begin{equation}
\label{Bedicek}
F^q {\it Hom}(M,N) := \{f \in {\it Hom}(M,N);\ f(F^pM) \subset
F^{p+q}N \mbox { for each $p$}\}.
\end{equation}
We believe that the notion of a filtered algebra, operad, etc.,~is
clear; we require that all structure operations preserve the
filtration.

If $M = (M,d_M)$ is a filtered chain complex, then~(\ref{Bedicek})
defines a filtration of the endomorphism
operad $\End_M$. A {\em
filtered algebra\/} over a filtered operad $\calP$ is a homomorphism
$A : \calP \to \End_M$ of filtered operads.
There is an evident colored analog of this notion.

\begin{example}
\label{Jojojoj}
{\rm
Let $\xb$ be a generator of degree $-1$ and let
\begin{equation}
\label{Ciperka}
\preDif := ({\calF(\xb)}, d),       
\end{equation}
with $d$ the `derivation' in the operadic sense defined by  
$d \xb := -\xb\, \xb$. 
The free operad $\calF(\xb)$ on $\xb$ is the same as the polynomial
ring ${\Bbb Z}[\xb]$. We
define the filtration
\[
F^p \preDif := \mbox { the subspace spanned by monomials in $\xb$ of
length $\geq p \geq 0$.}
\]
 
The differential $d$ clearly preserves the filtration, as well as does
the composition, so the operad $\preDif$ is filtered. 
Let $\Dif$ be the completion of $\preDif$ with respect to the above
filtration; of course, $\Dif$ coincides with the algebra of  power series
${\Bbb Z}[[\xb]]$.

Filtered $\Dif$-algebras $A : \Dif \to \End_M$ on $M = (M,d_M)$ 
correspond to perturbations $\td_M = d_M + \pa_M$ of the differential
$d_M$, the correspondence being given by 
$\pa_M := A(\xb)$. 
Indeed, $d \xb = - \xb\,\xb$ is mapped by
$A$ to $\pa_Md_M + d_M\pa_M = - \pa_M\pa_M$, which is the same as
$(d_M + \pa_M)^2=0$.
}
\end{example}

\begin{proposition}
\label{Jisasek}
The operad $\Dif$ of Example~\ref{Jojojoj} is acyclic, that is,
$H_*(\Dif) \cong {\frak 1}$, where 
${\frak 1}$ is the trivial operad.
\end{proposition}

The {\bf proof\/} of the above proposition is easy and we leave it as an
exercise. One feels that the proposition must be
`philosophically' true. Algebras over the trivial operad ${\frak 1}$
are just chain complexes with no additional structure, i.e.~with
only the structure given by the unperturbed
differential. The operad $\Dif$
describes perturbations of this differential, so it must be
homologically the same as ${\frak 1}$.

\section{Strong homotopy equivalences}
\label{Strong}

In Example~\ref{trochu} we introduced a colored operad $\Iso$
describing chain maps $\F: M \to N$, $\G: N\to M$
such that $\F\G = \idmala_N$ and $\F\G = \idmala_M$.

A general belief is that the homotopy analog of this
situation is given by a
quadruple $(\F,\G,\H,\L)$, where $\F : M \to N$ and $\G : N \to M$ are
degree zero chain maps that are homotopy inverses of each other, with
associated homotopies $\H$ and $\L$:
\begin{equation}
\label{boli_mne_zada}
\G\F - \idmala_M = d_M\H + \H d_M,\ \F\G - \idmala_N = d_N\L + \L d_N.
\end{equation}
Such a quadruple is clearly an algebra over the operad
\begin{equation}
\label{Klapacek}
\calRfake  := (\calF(f_0,\ g_0,\ f_1,\  g_1),d),
\end{equation}
where
\[
f_0 : \B \to \W,\
g_0 : \W \to \B,\
f_1 : \B \to \B \mbox { and }
g_1 : \W \to \W
\]
are generators with 
$\deg(f_0) = \deg(g_0) = 0$, $\deg(f_1) = \deg(g_1) = 1$, and the
differential $d$ is given by
\[
d(f_0) := 0,\
d(g_0) := 0,\
d(f_1) := g_0f_0 - 1_\B
\mbox { and }
d(g_1) := f_0g_0 - 1_\W.
\]
There is a dg operad
map $\alphafake : \calRfake \to \Iso$ given by
\[
\alphafake(f_0) := f,\
\alphafake(g_0) := g,\
\alphafake(f_1) := 0
\mbox { and }
\alphafake(g_1) := 0.
\]

The following fact which shows that $\calRfake$ {\em is
not\/} an acyclic resolution of the operad $\Iso$ is crucial.

\vskip 2mm
\noindent
{\bf Fact.}
{\em
The map $\alphafake$ is not a homology isomorphism. For instance,
$f_0f_1 - g_1g_0$ is a cycle in the kernel of $\alphafake$ that is not
homologous to zero.
}

A proper resolution of $\Iso$ was described in~\cite{markl:ha}.
It is a graded colored differential operad
\[
\Riso := (\calF(f_0,f_1,\ldots; g_0,g_1,\ldots), d),
\]
with generators of two types,
\begin{equation}
\label{Myska}
\def\arraystretch{1.2}
\begin{array}{rl}
\mbox{(i)}&\mbox{%
\hskip -2mm 
generators $\seq{f_n}{n \geq 0}$, $\deg(f_n)=n$,} 
\left\{
\begin{array}{l}
\mbox{$f_n : \B \to \W$ if $n$ is even,} 
\\
\mbox{$f_n : \B \to \B$ if $n$ is odd,}
\end{array}
\right.
\\
\mbox{(ii)}&\mbox{%
\hskip -2mm
generators
$\seq{g_n}{n \geq 0}$, $\deg(g_n)=n$,}
\left\{
\begin{array}{l} 
\mbox{$g_n : \W \to \B$ if $n$ is
even,}
\\
\mbox{$g_n : \W \to \W$ if $n$ is odd.}
\end{array}
\right.
\end{array}
\end{equation}
The differential $d$ is given by
\[
\def\arraystretch{1.4}
\begin{array}{ll}
d f_0\widedef 0, 
& 
d g_0 \widedef 0, 
\\
d f_1 \widedef g_0f_0 - 1,   
& 
d g_1 \widedef f_0g_0 - 1    
\end{array}
\]
and, on remaining generators, by the formula
\begin{eqnarray}
\nonumber
df_{2m}   &:=& \sum_{0 \leq i < m}(f_{2i}f_{2(m-i)-1}-
g_{2(m-i)-1}f_{2i}),\ m \geq 0,
\\
\label{pejsek_a_kocicka}
df_{2m+1} &:=& \sum_{0 \leq j \leq m} g_{2j}f_{2(m-j)} -
               \sum_{0 \leq j < m} f_{2j+1}f_{2(m-j)-1},\ m \geq 1,
\\
\nonumber 
dg_{2m}   &:=& \sum_{0 \leq i < m}(g_{2i}g_{2(m-i)-1}-
f_{2(m-i)-1}g_{2i}),\ m \geq 0,
\\
\nonumber 
dg_{2m+1} &:=& \sum_{0 \leq j \leq m} f_{2j}g_{2(m-j)} -
               \sum_{0 \leq j < m} g_{2j+1}g_{2(m-j)-1},\ m \geq 1,
\end{eqnarray}
see also~(\ref{8.2}).
The above formulas can be written in a compact form by
introducing
elements
\begin{equation}
\label{bojim_se_porad}
\def\arraystretch{1.5}
\begin{array}{ll}
\fb := f_0 + f_2 + f_4 + \cdots : \B \to \W, &
\hb := f_1 + f_3 + f_5 + \cdots : \B \to \B,
\\
\gb := g_0 + g_2 + g_4 + \cdots : \W \to \B, &
\lb := l_1 + l_3 + l_5 + \cdots : \W \to \W.
\end{array}
\end{equation} 
Then $\Riso = \calF(\fb, \gb, \hb, \lb)$ with 
the differential given by
\[
d\fb = \fb\hb - \lb\fb,\   d\hb = \gb\fb - \hb\hb -1_\B,\
d\gb = \gb\lb - \hb\gb \mbox { and }   d\lb = \fb\gb - \lb\lb -1_\W. 
\]
We will use this kind of abbreviation quite often, but we
shall always keep in mind that each formula of this type in fact represents
infinitely many formulas for homogeneous parts.
The operad $\Riso$ is `trivially' filtered, by
\[
F^p\Riso := \cases{\Riso}{for $p = 0$, and}{0}{for $p > 0$.}
\]
This filtration is obviously complete. Algebras over the operad $\Riso$
are the strong homotopy equivalences introduced in
Definition~\ref{Kiki}. 

The following theorem, formulated without proof in~\cite{markl:ha}, 
claims that
$\Riso$ is an acyclic resolution of the operad $\Iso$.

\begin{theorem}
\label{Kukacky}
The map $\alphaiso : \Riso \to \Iso$ defined by
\begin{equation}
\label{Krtek_s_lopatickou}
\alphaiso(f_0) := [f],\ \alphaiso(g_0) := [g],
\mbox { while }
\alphaiso(f_n) := 0,\ \alphaiso(g_n) = 0 \mbox { for } n \geq 1,
\end{equation}
is a map of differential graded colored operads that induces an
isomorphism of cohomology.
\end{theorem} 

\noindent
{\bf Proof.}
It is clear that $\alphaiso$ commutes with the differentials and that
it induces an isomorphism $H_0(\Riso,d) \cong \Iso$. It thus remains 
to prove that $\Riso$ is acyclic in positive dimensions.

The operad $\calF(f_0,f_1,\ldots; g_0,g_1,\ldots)$ 
is the free abelian group spanned by composable 
chains of generators. The length
of these chains induces
another grading, which we call the
homogeneity. The differential $d$ decomposes as $d = d_{-1} +
d_{+1}$, where $d_{i}$ raises the homogeneity by $i= \pm 1$.
Explicitly, $d_{+1}$ is given on generators by
\[
d_{+1}\fb \!=\! \fb\hb - \lb\fb,\   d_{+1}\hb \!=\! \gb\fb - \hb\hb,\
d_{+1}\gb \!=\! 
\gb\lb - \hb\gb,\   d_{+1}\lb \!=\! \fb\gb - \lb\lb, 
\]
while $d_{-1}$ is given by
\[
d_{-1}\fb = 0,\   d_{-1}\hb =  -1_\B,\
d_{-1}\gb = 0 \mbox { and }   d_{-1}\lb = -1_\W. 
\]

We claim that 
\begin{equation}
\label{reboot}
(\calF(f_0,f_1,\ldots; g_0,g_1,\ldots),d_{+1}) 
\mbox{ is acyclic in positive degrees.}
\end{equation}
We prove~(\ref{reboot}) by introducing a contracting homotopy
\[
\theta
: \calF(f_0,f_1,\ldots; g_0,g_1,\ldots) \to \calF(f_0,f_1,\ldots;
g_0,g_1,\ldots)
\] 
as follows. Let $z_1,z_2,\ldots$ denote generators of
$\calF(f_0,f_1,\ldots; g_0,g_1,\ldots)$, then let, for $m \geq 0$,
\[
R(z_1z_2) :=
\left\{
\begin{array}{ll}
f_{2m+2}, & \mbox { if $z_1z_2 = f_0f_{2m+1}$,}
\\
g_{2m+2}, & \mbox { if $z_1z_2 = g_0g_{2m+1}$,}
\\
g_{2m+1}, & \mbox { if $z_1z_2 = f_0g_{2m}$,}
\\
f_{2m+1}, & \mbox { if $z_1z_2 = g_0f_{2m}$,}
\\
0 & \mbox { if otherwise.}
\end{array}
\right.
\]
Then the contracting homotopy $\theta$ is defined by
\[
\theta(z_1z_2 \cdots z_t) := 
\cases{R(z_1z_2)z_3 \cdots z_t}{if $t \geq 2$, and}
       {0}{otherwise.}
\]
It is immediate to check that indeed 
$\theta d_{+1}(x) + d_{+1} \theta(x) = x$
whenever $x$  has positive degree, which proves~(\ref{reboot}).

Suppose that $x \in \calF(f_0,f_1,\ldots; g_0,g_1,\ldots)$ 
is a $d$-cycle of
positive degree and let
\[
x = x_1 + \cdots + x_N, \mbox { $x_j$ has homogeneity $j$, $1 \leq j
\leq N$, $N > 1$,}
\] 
be its decomposition into
homogeneity-homogeneous parts. Then clearly $d_{+1}(x_N) = 0$, thus,
by~(\ref{reboot}), there exists some $b_{N-1}$ of homogeneity $N-1$ such
that $x_N = d_{+1}(b_{N-1})$. Then $x - d(b_{N-1})$ is a $d$-cycle
homologous to $x$, whose decomposition contains no terms of
homogeneity $\geq N$. Repeating this process as many times as
necessary, we end up with some $x'$, homologous to $x$, of homogeneity
$1$, i.e.~linear in the generators. An immediate inspection shows that
there is no such nontrivial $x'$ of positive degree, therefore $x'=0$
which finishes the proof, since $x'$ is, by construction, homologous
to $x$.%
\qed

Observe that, in the course of the proof of Theorem~\ref{Kukacky}, we
proved the following interesting statement:

\begin{proposition}
\label{Kravicka_Bu}
The map 
\[
\alpha_{+1} :
\left( 
\calF(f_0,f_1,\ldots; g_0,g_1,\ldots), d_{+1}
\right) 
\to
\frac{\calF(f,g)}{(fg =0,\ gf = 0)}
\]
given by 
\[
\mbox {
$\alpha_{+1}(f_0) := [f]$, $\alpha_{+1}(g_0) := [g]$, while
$\alpha_{+1}(f_n)= 0$ and $\alpha_{+1}(g_n)= 0$ for $n \geq 1$,
}
\]
is a homology isomorphism.
\end{proposition}

In the rest of this section we study when a homotopy
equivalence~(\ref{boli_mne_zada}) extends to a strong homotopy one.
Observe first that~(\ref{boli_mne_zada}) induces a `restricted' action
$\Ares : \Riso \to \End_{M,N}$ by
\[
\Ares(f_0) = \F,\
\Ares(g_0) = \G,\
\Ares(f_1) = \H \mbox { and } \Ares(g_1) =\L.
\]
Related to these data are two obstruction cycles
\begin{equation}
\label{kapesnik}
\frakO_M := \F\H -\L\F \in \Hom_1(M,N)
\mbox { and } \frakO_N := \G\L - \H\G \in \Hom_1(N,M).
\end{equation}

\begin{theorem}
\label{snad_nejsem_nastydly}
The obstruction $[\frakO_M] \in H_1(\Hom(M,N))$ 
vanishes if and only if the obstruction $[\frakO_N]
\in H_1(\Hom(N,M))$ does. 

The restricted action $\Ares$ can be extended to a full action 
$A : \Riso \to \End_{M,N}$ if and only if one (and hence both) of the
above obstructions vanish.
\end{theorem}

\noindent
{\bf Proof.}
Let us denote by $(\FF n,d)$ the suboperad of $\calF(f_0,f_1,\ldots;
g_0,g_1,\ldots)$ generated by $\{f_n,\ g_n\}_{i < n}$, 
with the induced differential. It is clear from the
definition that $df_n, dg_n \in \FF n$ for any $n \geq 1$, thus it
makes sense to consider the homology 
classes $[df_n]_{n-1}$ and  $[dg_n]_{n-1}$ of
these elements in $H_{n-1}(\FF n,d)$. We claim that
\begin{eqnarray}
\label{volne}
[g_0] [df_{2m}]_{2m-1} + [df_{2m}]_{2m-1} [f_0]& =& 0 
\mbox { in $H_{2m-1}(\FF { 2m},d)$,
and}
\\
\label{modely}
[f_0] [df_{2m+1}]_{2m} + [dg_{2m+1}]_{2m} [f_0] 
&=& 0 \mbox { in $H_{2m}(\FF { 2m +1},d)$.}
\end{eqnarray}
The first equation follows from the inspection of the degree $2m+1$
part of
\begin{equation}
\label{Mysacek}
d\hb^2 = d(\fb\gb)
\end{equation}
which is
\[
(d\hb^2)_{2m+1} =
f_0 (df_{2m}) + (dg_{2m}) g_0 + d(\sum_\doubless{i+j=m}{i,j  \geq 1}
f_{2i}g_{2j});
\]
equation~(\ref{Mysacek}) can be verified directly.
Equation~(\ref{modely}) follows in the same manner from 
\[
d(\fb \hb - \lb \fb) = 0.
\]
which follows from $d^2=0$.
Observe that~(\ref{volne}) gives, for $m=1$,
\[
[g_0] [df_{2}]_{2} + [df_{2}]_{2} [f_0] = 0
\]
which is mapped by $\Ares : \calF(f_{<2}; g_{<2}) \to \End_{M,N}$ 
to
\[
[\G][\F\H - \L\F] + [\G\L - \H\G][\F] = 0 \mbox { in $H_1(\Hom(M,N))$,}
\]
which is of course 
\[
[\G][\frakO_M] + [\frakO_N][\F] = 0.
\]
This implies the first part of the statement, since multiplication
by the homology class of $f$ (resp.~of $g$) is an isomorphism, as
these maps are homotopy invertible.

Let us prove the second part of the theorem. One implication is clear
-- if the restricted action
$\Ares$ can be extended to a full one, then obviously both
obstructions must vanish.

Suppose that both obstructions vanish. Then the restricted
action can be clearly
extended to $f_2$ and $g_2$, i.e.~on $\FF 3$;
we denote this extended
action by $A_2$.

Let us suppose that we have extended $\Ares$ to some 
$A_{n-1} : \FF n \to
\End_{M,N}$, $n \geq 3$, and try to extend it to $f_n$ and $g_n$. We must
distinguish whether $n$ is even or odd; suppose first that $n = 2m$.
The extension clearly exists if and only if 
\begin{equation}
\label{Disk}
\def\arraystretch{1.2}
\begin{array}{rl}
\mbox{$[A_{2m-1}(df_{2m})] = 0$}&
\mbox {in $H_{2m-1}(\Hom(M,N))$, and}
\\
\mbox {$[A_{2m-1}(dg_{2m})] = 0$}& 
\mbox {in $H_{2m-1}(\Hom(N,M))$.}
\end{array}
\end{equation}
This, unfortunately, need not be true in general, but we can use the
following trick. Observe that if we change the definition of
$A_{2m-1}(f_{2m-1})$ by adding a cycle $\phi \in \Hom_{2m-1}(M,M)$ and 
$A_{2m-1}(g_{2m-1})$ by adding a cycle $\psi \in \Hom_{2m-1}(M,M)$, 
the extension $A_{2m-1}$ remains well defined. 
We show that by such a `recalibration,' we may always
achieve that the elements in~(\ref{Disk}) vanish.
Indeed, it follows from the definition of the differential, from
$A_{2m-1}(f_0) = \F$ and $A_{2m-1}(g_0) = \G$, that~(\ref{Disk}) changes to 
\[
[\F][\phi] +  [A_{2m-1}(df_{2m})] - [\psi][\F]= 0
\mbox {  and }
[\G][\psi] + [A_{2m-1}(dg_{2m})] - [\phi][\G]= 0.
\] 
This system can clearly be solved if and only if
\[
[\G][A_{2m-1}(df_{2m})] + [A_{2m-1}(dg_{2m})][\F] = 0,
\] 
which is the image of~(\ref{volne}) under $A_{2m-1}$. The case of odd $n$ is
discussed in the same manner, using~(\ref{modely}) instead
of~(\ref{volne}).%
\qed

In the light of Theorem~\ref{snad_nejsem_nastydly}, we will make no
distinction between $[{\frak o}_M]$ and $[{\frak o}_N]$ and denote both
obstructions by $[\frak o]$.
The following statement is a `chain-level' version of 
Theorem~\ref{snad_nejsem_nastydly}.

\begin{theorem}
\label{Marcelka}
The restricted action $\Ares$ can be extended to a 
full action $A : \Riso \to \End_{M,N}$ by putting $\Ares(f_n) = 0$ and
$\Ares(g_n) = 0$ for $n \geq 2$
if and only if the obstruction cycles~(\ref{kapesnik}) vanish and if
$\H\H = 0$ and $\L\L =0$.
\end{theorem}

The {\bf proof} is an easy exercise.
In~\cite{markl:ha} we formulated without proof the following
theorem.

\begin{theorem}
\label{3+1}
Let $(\F,\G,\H,\L)$ be a homotopy 
equivalence~(\ref{boli_mne_zada}). By changing either
$\H$ or $\L$ we may always achieve that the obstruction $[{\frak o}]$
vanishes, i.e.~that, by Theorem~\ref{snad_nejsem_nastydly}, 
the homotopy equivalence $(\F,\G,\H,\L)$
extends to a strong one. Examples of these changes are
\begin{eqnarray*}
(\F,\G,\H,\L) & \longmapsto & (\F,\G,\H-\G(\F\H-\L\F),\L),\mbox { or}
\\
(\F,\G,\H,\L) & \longmapsto & (\F,\G,\H,\L-\F(\G\L-\H\G)).
\end{eqnarray*}
\end{theorem}

\noindent
{\bf Proof.}
Let us show that the first substitution annihilates the obstructions. 
Denote for simplicity $\H' := \H -
\G(\F\H-\L\F)$. Then it can be verified directly that
\begin{eqnarray*}
{\frak o}_M(\F,\G,\H',\L) 
\hskip -2mm &=& \hskip -2mm 
\F\H - \L\F - \F\G(\F\H-\L\F) = d(-\L(\F\H-\L\F)),
\mbox { and }
\\
{\frak o}_N(\F,\G,\H',\L)  
\hskip -2mm &=&  \hskip -2mm
\G\L-\H\G + \G(\F\H-\L\F)g = d(\H^2 \G + \G \L^2 - \H\G\L), 
\end{eqnarray*}
therefore $[{\frak o}_M(\F,\G,\H',\L)] = [{\frak o}_N(\F,\G,\H',\L)] =0$.
The discussion of the second substitution is the same.%
\qed

\section{The retraction}
\label{Retraction}
\label{The_retraction}

Let us introduce a filtered colored operad $\tRiso$ describing
{\em perturbations\/} of strong homotopy equivalences. It is the
completion of the operad
$\pretRiso$
generated by two types of generators:
\begin{itemize}
\item[(i)]
generators $\{f_n\}_{n\geq 0}$ and $\{g_n\}_{n\geq 0}$ as
in~(\ref{Myska}) for an unperturbed strongly homotopy equivalence, and
\item[(ii)] generators for a perturbation, that is,
a generator $\xb$ for a perturbation of the
`black' differential, a generator $\yb$ for a perturbation of the
`white' differential, and generators $\fbar_n$ and $\gbar_n$ for
perturbations of $f_n$ resp.~$g_n$, $n \geq 0$.
\end{itemize}

For homogeneity of the notation, we will sometimes write $f^0_n$
(resp.~$g^0_n$) instead of $f_n$ (resp.~$g_n$) and $f^1_n$
(resp.~$g^1_n$) instead of $\fbar_n$  (resp.~$\gbar_n$). With these
conventions assumed,
\[
\pretRiso:= \left(
             \calF(\xb,\
                   \yb,\
                   \{f^s_n\}^{s=1,2}_{n \geq 0},\
                   \{g^s_n\}^{s=1,2}_{n \geq 0}
                  ),\ d
        \right)
\]
with $\deg(\xb) = \deg(\yb) = -1$ and $\deg(f_n^s) = \deg(g_n^s) =
n$. The differential $d$ will be defined later.
To define on $\pretRiso$ a filtration, we assign to each generator another
degree $\udeg$ by
\[
\udeg(\xb) = \udeg(\yb) = \udeg(\fbar_n) = \udeg(\gbar_n) = 1,\
\udeg(f_n) = \udeg(g_n)=0,\ n \geq 0.
\]
This assignment expresses the fact that overlined generators describe
perturbations. The $\udeg$-grading of generators induces, in the
standard way, a grading on $\pretRiso$ and we define
\[
F^p \pretRiso := \{z \in \pretRiso;\ \udeg(z) \geq p\},\ p \geq 0.
\]
Let us denote by $\tRiso$ the completion of $\pretRiso$. 
A typical element of $\tRiso$ is a formal sum $\sum_{i \geq 0}m_i$
with $m_i \in \pretRiso$ and $\udeg(m_i) = i$. 

The best way to describe the differential is to introduce a condensed
notation (compare~(\ref{bojim_se_porad})):
\begin{equation}
\label{Konicek}
\def\arraystretch{1.5}
\begin{array}{ll}
\fbtilde := \sum_{m\geq 0}f_{2m} + \fbar_{2m},
&
\gbtilde := \sum_{m\geq 0}  g_{2m} + \gbar_{2m},
\\
\hbtilde := \sum_{m\geq 0}f_{2m+1}+ \fbar_{2m+1},\
&
\lbtilde := \sum_{m\geq 0}  g_{2m+1} + \gbar_{2m+1}.
\end{array}
\end{equation}
The differential is given by
\begin{equation}
\label{Brumlalek}
\def\arraystretch{1.5}
\begin{array}{ll}
d\xb \!=\! -\xb\,\xb, & d\yb \!=\! - \yb\,\yb,
\\
d\fbtilde \!=\! \fbtilde(\xb + \hbtilde) \!-\! (\yb + \lbtilde)\fbtilde, & 
d\gbtilde \!=\! \gbtilde(\yb + \lbtilde) \!-\! (\xb + \hbtilde)\gbtilde,
\\
d\hbtilde \!=\! \!-\! (\hbtilde \xb \!+\! \xb\hbtilde) \!+\! 
\gbtilde\fbtilde \!-\! \hbtilde\hbtilde \!-\!1, &
d\lbtilde \!=\! \!-\! (\lbtilde \yb \!+\! \yb\lbtilde) \!+\! 
\fbtilde\gbtilde \!-\! \lbtilde\lbtilde \!-\!1.
\end{array}
\end{equation}
A moment's reflection shows that the differential operad $\tRiso$
really describes
perturbations of strongly homotopy equivalences. Expanding~(\ref{Brumlalek})
we get more explicit formulas for the differential:
\begin{eqnarray*}
df^1_{2m} &:=& \sum_{t=1,2}
      (f^{t}_{2m}\xb - \yb f^{t}_{2m}) + \sum_{t+r \geq 1}(
      \sum_{0 \leq i < m}(f^t_{2i}f^{r}_{2(m-i)-1}- 
                          g^t_{2(m-i)-1}f^{r}_{2i})), 
\\
df^1_{2m+1} &:=& \sum_{t=1,2}
     -(f^{t}_{2m+1}\xb + \xb f^{t}_{2m+1}) +
\\&&+
             \sum_{t+r \geq 1}(
             \sum_{0 \leq j \leq m} g^t_{2j}f^{r}_{2(m-j)} -
             \sum_{0 \leq j < m} f^t_{2j+1}f^{r}_{2(m-j)-1}),
\\
dg^1_{2m} &:=& \sum_{t=1,2}
      (g^{t}_{2m}\yb - \xb g^{t}_{2m}) + \sum_{t+r \geq 1}(
      \sum_{0 \leq i < m}(g^t_{2i}g^{r}_{2(m-i)-1}- 
                          f^t_{2(m-i)-1}g^{r}_{2i})), 
\\
dg^1_{2m+1} &:=& \sum_{t=1,2}
     -(g^{t}_{2m+1}\yb + \yb g^{t}_{2m+1}) +
\\&&+
             \sum_{t+r \geq 1}(
             \sum_{0 \leq j \leq m} f^t_{2j}g^{r}_{2(m-j)} -
             \sum_{0 \leq j < m} g^t_{2j+1}g^{r}_{2(m-j)-1}).
\\
\end{eqnarray*}
The action of $d$ on $f^0_n = f_n$ and $g^0_n = g_n$ is, of course,
the same as in~(\ref{pejsek_a_kocicka}). See also~\ref{8.3}.
The following theorem claims that $\tRiso$ is a resolution of the
operad $\Iso$ introduced in Example~\ref{trochu}.

\begin{theorem}
\label{mozna}
The map $\alpha : \tRiso \to \Iso$ given by
\begin{equation}
\label{Krtek_a_zrcatko}
\alpha(f^0_0) := [f] \mbox { and } \alpha(g^0_0) := [g],
\end{equation}
while $\alpha$ is zero on the remaining generators, is a map of
differential filtered operads that induces an isomorphism of cohomology.
\end{theorem}

\noindent 
{\bf Proof.}
It is immediate to see that $\alpha$ decomposes as $\alpha = \alphaiso
\circ \tilde{\alpha}$, with $\tilde{\alpha} : \tRiso \to \Riso$ 
given by $\tilde{\alpha}(f^0_n)= f_n$,
$\tilde{\alpha}(g^0_n) = g_n$ and $\tilde{\alpha}$ trivial on
remaining generators.

Since $\alphaiso$ is, by Theorem~\ref{Kukacky}, a homology isomorphism, it is
enough to show that $\tilde{\alpha}$ is also a homology isomorphism. 
This can be done by a spectral sequence argument which we
omit, since we will not need the theorem in our proofs. 

The philosophical meaning is that a perturbation cannot
introduce nontrivial homology classes.%
\qed

\vskip -3mm
Let us consider the free product
\[
\freeprod = (\calF(\xb,\seq {f_n}{n\geq 0}, \seq {g_n}{n\geq 0}),\ d)
\]
with the differential given by~(\ref{Ciperka}) and~(\ref{pejsek_a_kocicka}).
It is clear that the map $\iota :\freeprod \hookrightarrow \tRiso$
defined by 
\begin{equation}
\label{se}
\iota(f_n) := f^0_n,\
\iota(g_n) := g^0_n \mbox { and } \iota(\xb) := \xb,\ n \geq 0, 
\end{equation}
is an inclusion of differential filtered colored operads.
Let us formulate the main statement of this section.  

\begin{theorem}
\label{Motorka}
The operad $\freeprod$ is a retract of $\tRiso$, that is, 
there exists a map $r : \tRiso \to \freeprod$ of differential filtered
colored operads such that $r\iota = \id_{\freeprod}$. 
\end{theorem}

\noindent
{\bf Proof.}
We construct the retraction $r$ explicitly. Let us define, for each
odd $r \geq -1$, a
`kernel' $\kernel_r :  \B \to \B$, $\kernel_r \in\freeprod$, 
of degree $r$ by the formula
\[
\kernel_r := 
\sum_{t \geq 0} 
\xb f_{2m_1+1} \xb
\cdots \xb f_{2m_{t}+1} \xb
\]
where the summation runs over all 
$2(m_1+\cdots + m_{t}) -1 =r$, 
$m_1 \geq 0, \ldots, m_{t}\geq 0$. See~\ref{8.5} 
for some explicit formulas.
The retraction $r : \tRiso \to \freeprod$ is then given by the
following formulas:
\begin{equation}
\label{Mikinka}
\def\arraystretch{1.5}
\begin{array}{rclrcl}
r(\xb) & \hskip -2mm := \hskip -2mm& \xb, 
&  r(\yb)& \hskip -2mm := \hskip -2mm& f_0 \kernel_{-1} g_0,
\\
r(f_n)& \hskip -2mm := \hskip -2mm &f_n,& r(g_n) 
& \hskip -2mm := \hskip -2mm&  g_n, 
\\
r(\fbar_{2m})& \hskip -2mm 
:=\hskip -3mm \hskip -2mm& 
\displaystyle\sum_{a+b+c = m}\hskip -2mm  f_{2a}\kernel_{2b-1}f_{2c+1}, &
r(\gbar_{2m}) &\hskip -2mm 
:=\hskip -3mm \hskip -2mm &
\displaystyle\sum_{a+b+c = m}\hskip -2mm f_{2a+1} \kernel_{2b-1} g_{2b},
\\
r(\fbar_{2m+1}) &\hskip -2mm 
:=\hskip -3mm \hskip -2mm& 
\displaystyle\sum_{a+b+c = m} \hskip -2mm f_{2a+1}\kernel_{2b-1} f_{2c+1},&
r(\gbar_{2m+1})& \hskip -2mm 
:=\hskip -3mm \hskip -2mm& 
\displaystyle\sum_{a+b+c = m+1}\hskip -2mm f_{2a}\kernel_{2b-1}g_{2b},
\end{array}
\end{equation}
where $m,n \geq 0$ and $a,b,c$ are nonnegative integers.
In compact notation
\[
\kernelb := \sum_{q \geq 0} (\xb \hb)^q \xb
\] 
we can rewrite~(\ref{Mikinka}) as 
\[
\def\arraystretch{1.4}
\begin{array}{ll}
r(\fbtilde) = \fb(1 + \kernelb\hb), &
r(\gbtilde) = (1 + \hb\kernelb)\gb,
\\
r(\hbtilde) = \hb+\hb\kernelb\hb, &
r(\yb + \lbtilde) = \lb + \fb\kernelb\gb,
\end{array}
\]
see~(\ref{Konicek}) for the meaning of $\fbtilde, \gbtilde, \hbtilde$ and
$\lbtilde$. 
It is clear that $r$ defined above is a retraction. Let us prove that
it commutes with the differentials, that is
\begin{equation}
\label{Hopsalek}
dr = rd.
\end{equation}
It is, of course, enough to prove~(\ref{Hopsalek}) on generators
$\fbtilde,\ \gbtilde,\ \hbtilde$ and $\lbtilde$ of $\tRiso$. 
For $\fbtilde$ we have
\begin{eqnarray}
\label{son}
dr(\fbtilde) &=&
d(\fb(1 +\kernelb\hb)) = (\fb \hb - \lb\fb)(1+\kernelb\hb) +
\\
\nonumber 
&&+ \fb\kernelb(\gb\fb-
\hb\hb)\kernelb\hb - \fb\kernelb(\gb\fb-\hb\hb-1),
\end{eqnarray}
where we used the obvious relation
\[
d\kernelb = -\kernelb(\gb\fb - \hb\hb)\kernelb. 
\]
On the other hand, 
\begin{eqnarray}
\label{SON}
rd(\fbtilde) &=& r(\fbtilde(\xb + \hbtilde)-(\yb+\lbtilde)\fbtilde)=
\\
\nonumber 
&=&
\fb(1+\kernelb\hb)(\xb+\hb+\hb\kernelb\hb) - 
(\lb{} + \fb\kernelb\gb)\fb(1 + \kernelb\hb).
\end{eqnarray}
Comparing~(\ref{son}) to ~(\ref{SON}), using another obvious relation
\[
\xb + \kernelb \hb\xb = \kernelb,
\]
we indeed check that $rd(\fbtilde)=
dr(\fbtilde)$. Equation~(\ref{Hopsalek}) can be verified on
remaining generators by the same direct argument.%
\qed

\section{Proofs}
\label{Proofs}

The initial data of the Perturbation Problem  define an algebra 
$\Ain$ over the free product
\[
\Dif * \calRfake = (\calF(\xb,f_0,g_0,f_1,g_1),d),
\]
of the operad $\Dif$ of Example~\ref{Jojojoj} 
with the operad $\calRfake$
introduced in~(\ref{Klapacek}), $\Ain: \Dif * \calRfake \to \End_{M,N}$, by
\[
\Ain(\xb) := \pa_M := \td_M - d_M,\
\Ain(f_0) := \F,\ \Ain(g_0) := \G,\ \Ain(f_1) := \H,\
\Ain(g_1) := \L.
\]

We seek a solution of the PP encoded to an algebra 
$\Aout$ over the differential filtered
suboperad 
\[
\tcalRfake :=
(\calF(\xb,\yb,f_0,f_1,g_0,g_1,\fbar_0,\fbar_1,\gbar_0,\gbar_1), d)
\]
of the operad $\tRiso$ introduced in Section~\ref{Retraction} as 
\begin{eqnarray*}
&\td_N := \Aout(\yb),\
\tF \:= \Aout(\fbar_0) + \F,\
\tG \:= \Aout(\gbar_0) + \G, &
\\
&\tH := \Aout(\fbar_1) + \H \mbox { and }
\tL \:= \Aout(\gbar_1) + \L. &
\end{eqnarray*}

There is a natural
inclusion
$\iotafake: \Dif * \calRfake \hookrightarrow \tcalRfake$
given by
\[
\iotafake(\xb) := \xb,\
\iotafake(f_0) := f_0,\
\iotafake(g_0) := g_0,\
\iotafake(f_1) := f_1 \mbox { and }
\iotafake(g_1) := g_1.
\]
A `functorial' solution of the Perturbation Problem means to find
a retraction
\begin{equation}
\label{zase}
\rfake : \tcalRfake \to \Dif * \calRfake,\ \rfake\iotafake = 
\id_{\Dif * \calRfake}.
\end{equation}

\begin{theorem}
\label{Harvey}
There is no retraction $\rfake :  \tcalRfake \to \Dif * \calRfake$ 
as in~(\ref{zase}).
\end{theorem}

\noindent
{\bf Proof.}
The proof is a straightforward obstruction theory, but since the
non-existence of the retraction $\rfake$ motivated all this work, we
reproduce the proof here in its full length. All calculations below are
made modulo terms of filtration $\geq 2$, so we, in fact, work in the
associated graded operad.
The following equations must be satisfied 
(see~\ref{8.5}):
\begin{eqnarray}
\nonumber
d \rfake(\yb) \hskip -3mm &=& \hskip -3mm \rfake (d \yb) = 0,
\\
\label{1}
d \rfake(\fbar_0) \hskip -3mm &=& \hskip -3mm  \rfake(d \fbar_0) = f_0 \xb - \rfake(\yb) f_0,
\\
\label{2}
d \rfake(\gbar_0) \hskip -3mm &=& \hskip -3mm  \rfake(d \gbar_0) = g_0 \rfake(\yb) - \xb g_0, 
\\
\label{3}
d \rfake(\fbar_1) \hskip -3mm &=& \hskip -3mm  \rfake(d \fbar_1) =
                   -(\xb f_1 + f_1 \xb) + \rfake (\gbar_0)f_0 + g_0
                   \rfake(\gbar_0), \mbox { and}
\\
\label{4}
d \rfake(\gbar_1) \hskip -3mm &=& \hskip -3mm  \rfake(d \gbar_1) =
                   -(\rfake(\yb) g_1 + g_1 \rfake(\yb)) +
                   \rfake(\fbar_0)g_0 + f_0 \rfake(\gbar_0).
\end{eqnarray}
It follows from~(\ref{1}) and~(\ref{2}) that, for some $b$,
$\rfake(\yb) = f_0 \xb g_0 + db$
and that
\[
\rfake(\fbar_0) = f_0 \xb g_0 - bf_0 +c_1,\
\rfake(\gbar_0) = f_1 \xb g_0 + g_0b  + c_2,
\]
for some cycles $c_1,c_2$.
The right hand side of~(\ref{3}) then becomes
\[
(g_0f_0 - 1)\xb f_1 - f_1 \xb (1- g_0f_0) + c_2f_0 + g_0 c_1 
= d(f_1\xb f_1) + c_2 f_0 + g_0 c_1,
\]
while the right hand side of~(\ref{4}) becomes
\begin{eqnarray*}
\lefteqn{\hskip -1.8cm
f_0\xb(f_1g_0 - g_0g_1) + (f_0f_1 - g_1 f_0)\xb g - b f_0g_0 - db g_1
+ f_0g_0 b - g_1 db + c_1 g_0 + f_0 c_2 =
}
\\
&=& d(g_1b - bg_1) + f_0\{\xb(f_1g_0-g_0g_1) + c_2\} +
\{(f_0f_1-g_1f_0)\xb + c_1\}g_0. 
\end{eqnarray*}
{}From this we see that~(\ref{3})
and~(\ref{4}) can be solved in $\rfake(\fbar_1)$ and $\rfake(\gbar_1)$
if and only if 
\[
f_0\xb(f_1g_0-g_0g_1) + (f_0f_1-g_1f_0)\xb g_0
\]
is homologous to zero. It can be easily seen that this is not true.%
\qed

\noindent
{\bf Proof of the IPL.\/}
The initial data of the IPL can be organized into an action $E_{\rm
in} : \Dif *\Riso \to \End_{M,N}$. Then the action
\[
\Eout : \tRiso
\stackrel{r}{\longrightarrow} \freeprod \stackrel{E_{\rm in}}{\vlra}
\End_{M,N},
\]
where $r$ is the retraction of Theorem~\ref{Motorka}, clearly solves the
IPL.%
\qed

\section{A conceptual explanation}
\label{conceptual}

We believe in the existence of a model category (MC) structure on
the category of operads. Let us ignore in this conceptual section the
fact that the existence of this structure has been proved only for some
special cases~\cite{hinich:CA97} 
and certainly not for the category $\FiltOp$ of general filtered
colored operads over 
${\Bbb Z}$.

Our candidate for cofibrations in $\FiltOp$ are maps such that the
associated maps of graded operads are 
cofibrations in the sense of an obvious integral
version of~\cite[Definition~15]{markl:ha} (or something close to it).
Fibrations are then epimorphisms and weak equivalences are homology
isomorphisms. 

As argued in~\cite{markl:ha}, homotopy invariant algebras are those
over cofibrant operads. By Theorem~\ref{Kravicka_Bu}, 
$\Riso$ is a cofibrant
resolution of the operad $\Iso$, that is why strong homotopy
equivalences, as algebras over $\Riso$, are proper homotopy
versions of strict isomorphisms.

Let us show how the IPL follows from the properties of the MC
structure on $\FiltOp$.
The situation is summarized in the following diagram.

\begin{center}
{
\unitlength=1.2pt
\begin{picture}(230.00,75.00)(20.00,0.00)
\put(60.00,70.00){\makebox(0.00,0.00){$\id$}}
\put(0.00,30.00){\makebox(0.00,0.00){$\iota$}}
\put(60.00,10.00){\makebox(0.00,0.00){$\alpha$}}
\put(140.00,30.00){\makebox(0.00,0.00){$p$}}
\put(190.00,70.00){\makebox(0.00,0.00){$E_{{\rm in}}$}}
\put(130.00,0.00){\makebox(0.00,0.00){$\Iso$}}
\put(10.00,0.00){\makebox(0.00,0.00){$\tRiso$}}
\put(130.00,60.00){\makebox(0.00,0.00){$\Dif*\Riso$}}
\put(10.00,60.00){\makebox(0.00,0.00){$\Dif*\Riso$}}
\put(260.00,60.00){\makebox(0.00,0.00){$\End_{M,N}$}}

\multiput(16,3)(3,1.5){31}{\makebox(0,0){$\cdot$}}
\put(110.00,50.00){\vector(2,1){1}}
\put(70.00,40.00){\makebox(0.00,0.00){$r$}}

\put(10.00,50.00){\vector(0,-1){40.00}}
\put(130.00,50.00){\vector(0,-1){40.00}}

\put(30.00,0.00){\vector(1,0){80.00}}
\put(32.00,60.00){\vector(1,0){72.00}}
\put(160.00,60.00){\vector(1,0){70.00}}

\end{picture}}
\end{center}

In the above diagram,
$\alpha$ is the map from Theorem~\ref{mozna}, $\iota$ the
inclusion~(\ref{se}), 
$p := \alpha \circ \iota$ and the action 
$E_{\rm in}$ summarizes the input data of the IPL. The solution of the
IPL will then be given by $E_{\rm out} : = E_{\rm in} \circ r$.

The map is
$p$ clearly an epimorphism, hence a fibration.
It is also a weak equivalence, because, if we ignore the acyclic
(by Proposition~\ref{Jisasek})
factor $\Dif$, the map $p$ is exactly the map $\alphaiso$ of
Theorem~\ref{Kukacky}. The map
$\iota$ is a cofibration, thus the existence of $r$ follows from the 
axioms of a MC structure.

The above is, of course, just an explanation, not a proof, so we had,
in the proof of Theorem~\ref{Motorka},
to construct the retraction $r$ by other means.

\section{Appendix: Explicit Formulas}
\label{Appendix}

\begin{odstavec}
\label{8.1}
Explicit axioms for a SHE:
\begin{eqnarray*}
d_N F_0 - F_0 d_M  &=&  0,
\\
d_M G_0 - G_0 d_N  &=&  0,
\\
d_M H_1 + H_1 d_M  &=&  G_0 F_0 - \id_M, 
\\
d_N L_1 + L_1 d_N  &=&  F_0 G_0 - \id_N,
\\
d_N F_2 - F_2 d_M  &=&  F_0 H_1 - L_1 F_0,
\\
d_M G_2 - G_2 d_N  &=&  G_0 L_1 - H_1 G_0,
\\
d_M H_3 + H_3 d_M  &=&  G_0 F_2 - H_1 H_1 + G_2 F_0,
\\
d_N L_3 + L_3 d_N  &=&  F_0 G_2 - L_1 L_1 + F_2 G_0,
\\
d_N F_4 - F_4 d_M  &=&
F_0 H_3 - L_1 F_2 + F_2 H_1 - L_3 F_0,
\\
d_M G_4 - G_4 d_N &=& 
G_0 L_3 - H_1G_2  + G_2 L_1 - H_3 G_0,
\\
&\vdots&
\end{eqnarray*}
\end{odstavec}

\begin{odstavec}
\label{Alicek}
\label{8.2}
Formulas for the differential of $\Riso$:
\[
\def\arraystretch{1.4}
\begin{array}{rclrcl}
df_0 &=&0,&
dg_0 &=&0,
\\
df_1 &=& g_0f_0 -1,&
dg_1 &=& f_0g_0 -1,
\\
d f_2 &=& f_0f_1 - g_1f_0, 
&
d g_2 &=& g_0g_1 - f_1g_0, 
\\
d f_3 &=& g_0f_2 - f_1f_1 + g_2f_0,  
&
d g_3 &=& f_0g_2 - g_1g_1 + f_2g_0
\\
d f_4 &=& f_0f_3 - g_1f_2 + f_2f_1 
                      - g_3f_0, \hskip 3mm
&
d g_4 &=& g_0g_3 - f_1g_2 + g_2g_1 
                      - f_3g_0,
\\
&\vdots&&&\vdots&
\end{array}
\]
\end{odstavec}

\begin{odstavec}
\label{8.3}
Formulas for the differential of $\tRiso$: {\rm the action on $f_n$,
$g_n$ is, for $n \geq 0$, the same as in~\ref{Alicek}, and}
\begin{eqnarray*}
d \xb \hskip -2mm &=& 
\hskip -2mm - \xb\ \xb,\  d\yb = - \yb \ \yb,
\\
d\fbar_0 \hskip -2mm &=& 
\hskip -2mm f_0 \xb - \yb f_0 + \fbar_0 \xb - \yb \fbar_0,
\\
d\gbar_0 \hskip -2mm &=& 
\hskip -2mm g_0 \yb - \xb g_0 + \gbar_0 \yb - \xb \hskip .5mm \gbar_0,
\\
d\fbar_1 \hskip -2mm &=& 
\hskip -2mm - (\xb f_1 + f_1 \xb) + (\gbar_0 f_0 + g_0 \fbar_0) -
           (\xb \fbar_1 + \fbar_1 \xb) + \gbar_0\fbar_0,
\\
d\gbar_1 \hskip -2mm &=& 
\hskip -2mm - (\yb g_1 + g_1 \yb) + (\fbar_0 g_0 + f_0 \gbar_0) -
           (\yb\hskip .5mm \gbar_1 + \gbar_1 \yb) + \fbar_0\gbar_0,
\\
d\fbar_2 \hskip -2mm &=& 
\hskip -2mm (f_2 \xb - \yb f_2) + (\fbar_0 f_1 + f_0 \fbar_1) - 
           (\gbar_1 f_0 + g_1 \fbar_0) + (\fbar_2 \xb - \yb \fbar_2) 
           + (\fbar_1 \fbar_0 - \gbar_1 \gbar_0),
\\
d\gbar_2 \hskip -2mm &=& 
\hskip -2mm (g_2 \yb - \xb g_2) + (\gbar_0 g_1 + g_0 \gbar_1) - 
           (\fbar_1 g_0 + f_1 \gbar_0) + (\gbar_2 \yb - \xb \gbar_2) 
           + (\gbar_1 \gbar_0 - \fbar_1 \fbar_0),
\\
&\vdots&
\end{eqnarray*}
\end{odstavec}

\begin{odstavec}
\label{8.4}
Formulas for the kernel $\kernel_n$:
\begin{eqnarray*}
\kernel_{-1}\hskip -2.5mm \hskip -1.3mm &=& 
\hskip -2mm \xb + \xb f_1 \xb + \xb f_1 \xb f_1 \xb + 
                 \xb f_1 \xb f_1 \xb f_1 \xb + \cdots,
\\
\kernel_{1} \hskip -1.3mm &=& 
\hskip -2mm \xb f_3 \xb + \xb f_1 \xb f_3 \xb + \xb f_3 \xb  f_1 \xb +
                \xb f_1 \xb f_1 \xb f_3 \xb  +
                \xb f_1 \xb f_3 \xb  f_1 \xb + \xb f_3 \xb  f_1 \xb
                 f_1 \xb + \cdots,
\\
\kernel_3  \hskip -1.3mm &=& 
\hskip -2mm \xb f_5 \xb + \xb f_5 \xb f_1 \xb + 
                \xb f_3 \xb f_3 \xb + \xb f_1 \xb f_5 \xb + 
                \xb f_1 \xb f_3 \xb f_3 \xb + \xb f_3 \xb f_1 \xb f_3
                 \xb + 
\\          && \hskip 5mm +
                \xb f_3 \xb f_3 \xb f_1 \xb + \xb f_1 \xb f_1 \xb f_5
                \xb  + \xb f_1 \xb f_5 \xb f_1 \xb +
                \xb f_5 \xb f_1 \xb f_1 \xb  + \cdots,
\\
&\vdots&
\end{eqnarray*}
\end{odstavec}

\begin{odstavec}
\label{8.5}
Formulas for the retraction $r: \tRiso \to \freeprod$:
\begin{eqnarray*}
r(\yb) &:=& f_0 \kernel_{-1} g_0
\\
r(\fbar_0) &:=& f_0 \kernel_{-1} f_1
\\
r(\gbar_0) &:=& f_1 \kernel_{-1} g_0
\\
r(\fbar_1) &:=& f_1 \kernel_{-1} f_1
\\
r(\gbar_1) &:=& f_0 \kernel_{-1}g_2 + f_2 \kernel_{-1}g_0 + f_0 \kernel_1 g_0
\\
r(\fbar_2) &:=& f_0 \kernel_{-1}f_3 + f_2 \kernel_{-1} f_1 + f_0 \kernel_{1} f_1
\\
r(\gbar_2) &:=& f_1\kernel_{-1}g_2 + f_3 \kernel_{-1} g_0 + f_1 \kernel_{1} g_0
\\
r(\fbar_3) &:=& 
           f_3 \kernel_{-1} f_1 + f_1 \kernel_{-1} f_3 + f_1 \kernel_{1} f_1
\\
r(\gbar_3) &:=& f_4 
            \kernel_{-1} g_0 + f_2 \kernel_{-1} g_2 + f_0 \kernel_{-1} g_4 +
           f_2 \kernel_{1} g_0 + f_0 \kernel_{1} g_2 + f_0 \kernel_{3}
           g_0, 
\\
&\vdots&
\end{eqnarray*} 
\end{odstavec}

\begin{odstavec}
\label{8.6}
Formulas for the solution of the IPL:
\begin{eqnarray*}
{\tilde d}_N \hskip -2mm &=& \hskip -2mm
d_N + \F_0(\pa_M + \pa_M \H_0 \pa_M + \pa_M \H_0 \pa_M \H_0 \pa_M + 
        \pa_M \H_0 \pa_M \H_0 \pa_M \H_0 \pa_M+\cdots )\G_0,
\\
\tilde \F_0 \hskip -2mm &=& \hskip -2mm \F_0 +\F_0(\pa_M + \pa_M \H_0 
       \pa_M + \pa_M \H_0 \pa_M \H_0 \pa_M + 
        \pa_M \H_0 \pa_M \H_0 \pa_M \H_0 \pa_M+ \cdots )\H_0,
\\
\tilde \G_0 \hskip -2mm &=& \hskip -2mm \G_0 +
               \H_0(\pa_M + \pa_M \H_0 \pa_M + \pa_M \H_0 \pa_M \H_0 \pa_M + 
        \pa_M \H_0 \pa_M \H_0 \pa_M \H_0 \pa_M+ \cdots )\G_0,
\\
\tilde \H_0 \hskip -2mm &=& \hskip -2mm \H_0 +\H_0(\pa_M + 
            \pa_M \H_0 \pa_M + \pa_M \H_0 \pa_M \H_0 \pa_M + 
        \pa_M \H_0 \pa_M \H_0 \pa_M \H_0 \pa_M+ \cdots )\H_0,
\\
\tilde \L_0 \hskip -2mm &=& \hskip -2mm \L_0 +\F_0(\pa_M + \pa_M \H_0 \pa_M + 
         \pa_M \H_0 \pa_M \H_0 \pa_M + 
        \pa_M \H_0 \pa_M \H_0 \pa_M \H_0 \pa_M+ \cdots )\G_2
\\
&& \hskip 3mm + \F_2(\pa_M + \pa_M \H_0 \pa_M + \pa_M \H_0 \pa_M \H_0 \pa_M + 
        \pa_M \H_0 \pa_M \H_0 \pa_M \H_0 \pa_M+ \cdots )\G_0
\\
&& \hskip 3mm+ \F_0(\pa_M \F_3 \pa_M + \pa_M \H_0 \pa_M \F_3 \pa_M + 
            \pa_M \F_3 \pa_M \H_0 \pa_M + \cdots)\G_0, 
\\
&\vdots&
\end{eqnarray*}
{\rm
In the above formulas, $\pa_M = \td_M - d_M$.}
\end{odstavec}

\begin{center}
\small \bf ACKNOWLEDGEMENTS
\end{center}
I would like to express my thanks to Jim Stasheff
and Johannes Huebschmann for reading the manuscript and many useful
comments, and apologize to the latter for not
implementing all his suggestions. 


\end{document}